\documentclass[12pt,leqno]{article}
\usepackage{amsmath,amssymb,amscd}
\numberwithin{equation}{section}

\newtheorem{proposition}{Proposition}[section]

\newtheorem{theorem}{Theorem}[section]
\newtheorem{corollary}{Corollary}[section]
\newtheorem{definition}{Definition}[section]
\newcommand\qed{{\unskip\nobreak\hfil\penalty50\hskip2em\vadjust{}
    \nobreak\hfil$\Box$\parfillskip=0pt\finalhyphendemerits=0\par}}

\renewcommand{\c}{\circ}
\newcommand{\x}{\times}
\renewcommand{\a}{\alpha}
\renewcommand{\b}{\beta}
\renewcommand{\d}{\delta}
\newcommand{\D}{\Delta}

\newcommand{\g}{\gamma}
\newcommand{\G}{\Gamma}
\renewcommand{\l}{\lambda}
\renewcommand{\L}{\Lambda}

\renewcommand{\th}{\theta}
\renewcommand{\O}{\Omega}
\renewcommand{\o}{\omega}

\renewcommand{\i}{\infty}

\newcommand{\bZ}{{\mathbb Z}}
\newcommand{\bE}{{\mathbb E}}
\newcommand{\bN}{{\mathbb N}}
\newcommand{\bP}{{\mathbb P}}
\newcommand{\bR}{{\mathbb R}}

\newcommand{\cA}{{\mathcal A}}

\newcommand{\cP}{{\mathcal P}}
\newcommand{\cF}{{\mathcal F}}

\newcommand{\cL}{{\mathcal L}}

\newcommand{\cC}{{\mathcal C}}
\newcommand{\cX}{{\mathcal X}}

\begin{document}

\title{Stochastically Symplectic Maps and Their Applications to Navier-Stokes Equation}
\author{Fraydoun Rezakhanlou\thanks{This work is supported in part by NSF Grant DMS-1106526.} \\
UC Berkeley \\
Department of Mathematics \\
Berkeley, CA\ \ 94720-3840}

\maketitle

\begin{abstract}
 Poincar\'e's invariance principle for Hamiltonian flows implies Kelvin's principle for solution to Incompressible Euler Equation. Iyer-Constantin Circulation Theorem offers a stochastic analog of Kelvin's principle for Navier-Stokes Equation. Weakly symplectic diffusions are defined to produce stochastically symplectic flows in a systematic way. With the aid of symplectic diffusions, we produce a family of martigales associated with solutions to Navier-Stokes Equation that in turn can be used to prove Iyer-Constantin Circulation Theorem. 
We also review some basic facts in symplectic and contact geometry
 and their applications to Euler Equation. 
\end{abstract}

\section{Introduction}
\label{sec1}

Hamiltonian systems appear in conservative problems of mechanics 
governing the motion of particles in fluid.  Such a mechanical system  is modeled by a
 Hamiltonian function $H(x,t)$ where $x = (q,p)\in {\mathbb R}^d \times {\mathbb R}^d$, $q = (q_1,\dots,q_d)$, $p = (p_1,\dots,p_d)$ denote the positions and the momenta of particles.  
The Hamiltonian's equations of motion are
\begin{equation}\label{eq1.1}
{\dot q} = H_p(q,p,t),\ {\dot p} = -H_q(q,p,t) 
\end{equation}
which is of the form
\begin{equation}\label{eq1.2}
{\dot x} = J\nabla_xH(x,t),\ \ \ J = \begin{bmatrix}
0 & I_d \\
-I_d & 0
\end{bmatrix} 
\end{equation}
where $I_d$ denotes the $d \times d$ identity matrix. 
It was known to Poincar\'e  that if $\phi_t$ is the flow of the ODE \eqref{eq1.2} and $\gamma$ is a closed curve, then
\begin{equation}\label{eq1.3}
\frac {d}{dt} \int_{\phi_t(\gamma)} \bar \l = 0,
\end{equation}
where $\bar\l:=p \cdot dq$.
 We may use Stokes' theorem to rewrite \eqref{eq1.3} as
\begin{equation}\label{eq1.4}
\frac {d}{dt} \int_{\phi_t(\Gamma)} {d\bar \l}   = 0 
\end{equation}
for every two-dimensional surface $\Gamma$.  In words, the $2$-form 
\[
{\bar \omega}:=\sum_{i=1}^d dp_i \wedge dq_i,
\]
 is invariant under the Hamiltonian flow $\phi_t$. Equivalently,
\begin{equation}\label{eq1.5}
\phi_t^*\bar \omega=\bar\o.
\end{equation}

 A Hamiltonian system \eqref{eq1.2} simplifies if we can find a function $u(q,t)$ such that $p(t) = u(q(t),t)$. 
 If such a function $u$ exists, then $q(t)$ solves
\begin{equation}\label{eq1.6}
\frac {dq}{dt} = H_p(q,u(q,t),t). 
\end{equation}
The equation for the time evolution of $p$ gives us an equation for the evolution of the velocity function $u$; since
\begin{align*}
{\dot p} &= (Du){\dot q} + u_t = (Du)  H_p(q,u,t) + u_t, \\
{\dot p} &= -H_q(q,u,t),
\end{align*}
the function $u(q,t)$ must solve,
\begin{equation}\label{eq1.7}
u_t +( Du)  H_p(q,u,t) + H_q(q,u,t) = 0. 
\end{equation}
For example, if $H(q,p,t) = \frac {1}{2} |p|^2 + P(q,t)$, then \eqref{eq1.7} becomes
\begin{equation}\label{eq1.8}
u_t + (Du)u + \nabla P(q,t) = 0,
\end{equation}
and the equation \eqref{eq1.6} simplifies to
\begin{equation}\label{eq1.9}
\frac {dq}{dt} = u(q,t) .
\end{equation}
 Here and below we write $Du$ and $\nabla P$ for the $q$-derivatives of the vector field $u$ and 
the scalar-valued function $P$ respectively. If the flow of \eqref{eq1.11} is denoted by $Q_t$, then $\phi_t(q,u(q,0)) = (Q_t(q),u(Q_t(q),t))$.  
Now \eqref{eq1.3} means that for any closed $q$-curve $\eta$, 
\begin{equation}\label{eq1.10}
\frac {d}{dt} \int_{Q_t(\eta)} u(q,t) \cdot dq =\frac {d}{dt} \int_{\eta}(DQ_t)^*\  u\circ Q_t(q,t) \cdot dq= 0,
\end{equation}
or equivalently
\begin{equation}\label{eq1.11}
d (Q_t^*\a_t)=d\a_0,
\end{equation}
where $\a_t=u(q,t)\cdot dq$.
This is the celebrated {\em Kelvin's circulation theorem}. In summary Poincar\'e's invariance principle 
\eqref{eq1.3}
implies Kelvin's principle for Euler Equation. (Note that the incompressibility condition $\nabla\cdot u=0$ is not
needed for \eqref{eq1.10}.) 

We may rewrite \eqref{eq1.11} as
\begin{equation}\label{eq1.12}
 Q_t^*(d\a_t)=d\a_0,
\end{equation}
and this is equivalent to Euler equation ( the equation \eqref{eq1.8} with the incompressibility condition $\nabla\cdot u=0$). Moreover, when $d=3$, \eqref{eq1.12} can be written as 

\begin{equation}\label{eq1.13}
\xi^t\circ Q_t=(DQ_t)\xi^0,\ \ \ {\text{ or }}\ \ \ \xi^t=\left((DQ_t)\xi^0\right)\circ Q_t^{-1},
\end{equation}
where $\xi^t(\cdot)=\nabla\times u(\cdot,t)$. The equation \eqref{eq1.13} is known as 
{\it Weber's formulation} of Euler Equation and is 
 equivalent to the vorticity equation by differentiating both sides with respect to $t$:
\begin{equation}\label{eq1.14}
\xi_t+(D\xi)u=(Du)\xi.
\end{equation}

 Constantin and Iyer [CI] discovered a circulation invariance principle 
for Navier-Stokes equation that is formulated in
 terms of a diffusion associated with the velocity field. Given a solution $u$ to the Navier-Stokes  equation 
\begin{equation}\label{eq1.15}
u_t +(Du)u + \nabla P(q,t) =  \nu   \D u , \ \ \ \nabla\cdot u=0,
\end{equation}
let us write $Q_t$ for the (stochastic) flow of the SDE
\begin{equation}\label{eq1.16}
dq=u(q,t)\ dt+\sqrt{2\nu }\ dW,
\end{equation}
with $W$ denoting the standard Brownian motion. If  we write $A=Q^{-1}$ and $\xi^t=\nabla\x u(\cdot,t)$,
and assume that $d=3$, then
Constantin and Iyer's circulation formula reads as
\begin{equation}\label{eq1.17}
\xi^t=\bE\left( (DQ_t)\xi^0\right)\circ A_t,
\end{equation}
where $\bE$ denotes the expected value.

We are now ready to state the first result of this article.
(To avoid a confusion between stochastic differential and exterior derivative, we use a hat for the latter.)

\begin{theorem}\label{th1.1} Write $\a_t=u(q,t)\cdot dq$ with $u$ a classical solution of \eqref{eq1.14}
and given $T>0$, set $B_t=Q_{T-t}\circ Q_T^{-1}$.
\begin{itemize}
\item(i)
Then the process $\b_t=B_t^*\hat d\a_{T-t}, \  t\in[0,T]$ is a 2-form valued martingale.
When $d=3$, this is equivalent to saying that the process
\[
M_t=\left(\left(DB_t^{-1}\right)\xi^{T-t}\right)\circ B_t,\ \ \ t\in[0,T],
\]
is a martingale.
\item (ii) Given a surface $\Theta$, the quadratic variation of the martingale
$\b_t({\Theta})=\int_\Theta \b_t$ is given by 
\[
\int_0^t\sum_{i=1}^d\left[\int_\Theta B_s^*\zeta_i^{T-s}\right]^2ds,
\]
 where
\[
\zeta_i^\th=\sum_{j,k=1}^d u^k_{q_iq_j}(\cdot,\th)\ dq_j\wedge dq_k,
\]
or equivalently, $\zeta_i^\th(v_1, v_2)=\cC(u_{q_i}(\cdot,\th))v_1\cdot v_2,$
with $\cC(w)=Dw-(Dw)^*$.
\item (iii) For $\Theta$ as in (ii), we have the bound
\[
\bE\int_0^T\sum_{i=1}^d\left[\int_\Theta B_s^*\zeta_i^{T-s}\right]^2ds\le 
\bE\left[\int_\Theta A_T^*\hat d\a_{0}\right]^2.
\]
\end{itemize}
\end{theorem}

\noindent
{\bf{Remark 1.1}}\begin{itemize}
\item (i) In a subsequent paper, we will show how
 Theorem 1.1 can be extended to certain weak solutions. To make sense of martingales $\b_t$ and $M_t$, 
we need to make sure that  $DQ_t$ exists weakly and belongs to suitable $L^r$ spaces. As it turns out, 
a natural condition to guarantee $DQ_t\in L^r$ for all $r\in[1,\i)$ is  
\[
\int_0^T\left[\int_{\bR^d}|u(x,t)|^p\ dx\right]^{q/p}\ dt<\i,
\]
 for some $p,q\ge 1$ such that $d/p+2/q\le 1.$
\item (ii) Our result takes a simpler form if $u$ is a solution to backward Navier-Stokes Equation. For such 
$u$, we simply have that $\b_t=Q_t^*\hat d\a_t$ is a martingale. When $d=3$, we deduce that 
$M_t=\big((DA_t)\xi^t\big)\circ Q_t$ is a martingale.
\end{itemize}
\qed

The organization of the paper is as follows: 
\begin{itemize}
\item In Section 2 we discuss Weber's formulation of Euler Equation and show how \eqref{eq1.5} implies \eqref{eq1.13}. 
We also discuss two fundamental results in Symplectic Geometry that are related to the so-called {\em{Clebsch}} variables.
\item In Section 3 we address some geometric questions for stochastic flows of general diffusions and
 study {\em {symplectic}} .
\item In Section 4 we use symplectic diffusions to establish Theorem 1.1.
\item In Section 5 we discuss {\em{contact}} diffusions.
\end{itemize}

\section{Euler Equation}
\label{sec2}
 \setcounter{equation}{0}

In this section we review some basic facts in differential geometry and their applications to Euler Equation.  Even though most of the discussion of this section is either well-known or part of folklore, a reader may find our discussion useful as we use similar ideas to prove Theorem 1.1.
We also use this section as an excuse to demonstrate/advertise  the potential use of symplectic/contact geometric
ideas in fluid mechanics.

We start with giving the elementary proof of \eqref{eq1.4}:
By Cartan's formula
\begin{equation}\label{eq2.1}
\frac {d}{dt}\phi^*_t\bar\l=\phi^*_t\cL_{Z_H}\bar\l=\phi^*_t dK=d(K\circ \phi_t),
\end{equation}
where  $\cL_Z$ denotes the Lie derivative with respect to the vector field $Z$,  $Z_H=J\nabla_xH$ for $H(q,p,t)=|p|^2/2+P(q,t)$, and
and  
\begin{equation}\label{eq2.2}
K(q,p,t)=p\cdot H_p(q,p,t)-H(q,p,t)=\frac 12|p|^2-P(q,t).
\end{equation}
If we integrate both sides of  \eqref{eq2.1} over an arbitrary (non-closed) curve of the form 
$(\eta,u(\eta,t))$, or equivalently restrict the form $\bar\l$ to the graph of the function $u$, then we obtain
\begin{equation}\label{eq2.3}
\frac {d}{dt}\left[(DQ_t)^*\  u\circ Q_t\right]=\nabla (L\circ Q_t),
\end{equation}
where $L(q,t)=K(q,u(q,t),t)=|u(q,t)|^2/2-P(q,t)$. 
Here by $A^*$ we mean the transpose of the matrix $A$. Recall $A_t=Q_t^{-1}$, so that 
\[
\left(DQ_t\right)^{-1}=DA_t\circ Q_t.
\]
As a consequence of \eqref{eq2.3} we have
\[
u(\cdot,t)= (DA_t)^*\ u^0\circ A_t+\nabla (R\circ A_t),\]
for $R=\int_0^t L\circ Q_s \ ds$. As a result,
\begin{equation}\label{eq2.4}
 u(\cdot,t)=\cP\left[ (DA_t)^*\ u^0\circ A_t\right],
\end{equation}
where $u^0$ is the initial data and 
$\cP$ denotes the Leray-Hodge projection onto the space of divergence-free vector fields.
The formula \eqref{eq2.4} is Weber's formulation and is equivalent to Euler's equation.

So far we have shown that the Kelvin's principle \eqref{eq1.3} 
is equivalent to the Weber's formulation of Euler equation. 
If we use \eqref{eq1.5} instead, we obtain a new equivalent 
formulation of Euler equation, namely the vorticity equation \eqref{eq1.12} or \eqref{eq1.13}. Recall
\[
\bar\o(v_1,v_2)=Jv_1\cdot v_2.
\]
If we choose $v_1$ and $v_2$ to be tangent to the graph of $u$, i.e. $v_i=(w_i,Du(q,t)w_i)$ for $i=1,2$, then
\[
\bar\o(v_1,v_2)=\cC(u)w_1\cdot w_2,
\]
where $\cC(u)=Du-(Du)^*$.  Hence \eqref{eq1.12} really means
\begin{equation}\label{eq2.5}
\cC(u(\cdot,t)\circ Q_t)(DQ_t)w_1\cdot(DQ_t)w_2=\cC(u(\cdot,0))w_1\cdot w_2.
\end{equation}

Let us assume now that $d=3$ so that, $\cC(u)w=\xi\x w$, where $\xi=\nabla\x u$ denotes the vorticity. 
Hence
\[
\bar\o(v_1,v_2)=(\xi\x w_1)\cdot w_2=:[\xi,w_1,w_2].
\]
We note that the right-hand side is the volume form evaluated at the triple $(\xi,w_1,w_2)$.  Now the invariance  
\eqref{eq2.5} becomes
\begin{equation}\label{eq2.6}
[\xi^t\circ Q_t,(DQ_t)w_1,(DQ_t)w_2]=[\xi^0,w_1,w_2],
\end{equation}
where we have written $\xi^t$ for $\xi(\cdot,t)$. Since $u$ is divergence-free, the flow $Q_t$ is volume preserving. 
As a result, 
\[
[\xi^0,w_1,w_2]=[(DQ_t)\xi^0 ,(DQ_t)w_1, (DQ_t)w_2].
\]
From this and \eqref{eq2.6} we deduce
\[
[\xi^t\circ Q_t,(DQ_t)w_1,(DQ_t)w_2]=[(DQ_t)\xi^0 ,(DQ_t)w_1,(DQ_t)w_2].
\]
Since $w_1$ and $w_2$ are arbitrary, we conclude that \eqref{eq1.13} is true.

\noindent
{\bf{Definition 2.1}}
 \begin{itemize}
\item(i) A closed $2$-form $\o$ is symplectic if it is nondegenerate. We say that symplectic forms
 $\o^1$ and $\o^2$ are {\em{isomorphic}} if there exists a diffeomorphism $\Psi$ such that $\Psi^*\o^1=\o^2$.
\item(ii) A $1$-form $\a$ is contact if $l_x=\{v:d\a(x;v,w)=0{\text{ for every }}w\}$ is a line and for every 
$v\in l_x$, we have that $\a(x;v)\neq 0$. 
We say that contact forms
 $\a^1$ and $\a^2$ are {\em{isomorphic}} if there exists a diffeomorphism $\Psi$ such that $\Psi^*\a^1=\a^2$.
We say that contact forms
 $\a^1$ and $\a^2$ are {\em{conformally isomorphic}} if there exist a diffeomorphism $\Psi$ and 
a scaler-valued continuous function $f>0$ such that $\Psi^*\a^1=f\a^2$.
\item(iii) A solution $u$ of Euler equation is called {\it{symplectic}} if $\o_0=d\a_0$ is symplectic.
\item(iv)  A solution $u$ of Euler equation is {\it{contact}} if there exists a scalar-valued $C^1$
 function $f_0$ such that $\a_0+df_0$ is contact. (Recall $\a_t=u(\cdot,t)\cdot dx$.)
\end{itemize}
\qed

\noindent
{\bf{Remark 2.1}}
 \begin{itemize}
\item(i) As it is well-known, the degeneracy of a 2-form can only happen when the dimension $d$ is even. 
Recall $\a_t=u(\cdot,t)\cdot dx$. If $u$ is a symplectic solution, then $\o_t=d\a_t$ is symplectic for all $t$
because by \eqref{eq1.12}, the form $\o_t$ is isomorphic to $\o_0$.
\item(ii)  When $u$ is a contact solution of Euler equation, then $\tilde \a_t=Q_t^*\a_0+df_t$ is contact for all $t$ where $f_t=f_0\circ Q_t$. In general $\tilde \a_t\neq \a_t$. However, by equation \eqref{eq1.12},  we have $d\a_t=d\tilde\a_t$. Hence there exists a scalar-valued function $g_t$ such that $\a_t+dg_t=\tilde \a_t$ is contact.
\end{itemize}
\qed

We continue with some general properties of symplectic and contact solutions of Euler Equation.

 As for symplectic solutions, assume that the dimension $d=2k$ is even and write
 \[
(q_1,\dots,q_d)=(x_1,y_1,\dots,x_k,y_k).
\] 
A classical theorem of Darboux asserts that all symplectic forms are isomorphic to the standard form
 $\bar\o=d\bar\l=\sum_{i=1}^kdy_i\wedge dx_i$. A natural question is whether such an isomorphism exists globally. 

\noindent
{\bf{Definition 2.2}}
Let $u$ be a symplectic solution of Euler Equation. We say that {\it Clebsch} variables exist for $u$
in the interval $[0,T]$, if we can find $C^1$ functions
$$X_1,\dots,X_k,Y_1,\dots,Y_k:\bR^d\times [0,T]\to\bR,\ \ \ F:\bR^d\times [0,T]\to\bR$$
such that $\Psi_t=(X_1,Y_1,\dots,X_k,Y_k)(\cdot,t)$ is a diffeomorphism, and 
$$
u(x,t)=\left(\sum_{i=1}^k Y_i\cdot \nabla X_i\right)(x,t)+\nabla F(x,t),
$$
for every $t\in[0.T]$. 
Alternatively, we may write $\a_t=\Psi_t^*\bar \l+dF$ or $d\a_t=\Psi_t^*\bar \o$.
\qed

\begin{proposition}\label{prop2.1} Let $u$ be a symplectic solution to Euler Equation. 
\begin{itemize}
\item(i) If Clebsch variables exist for $t=0$, then they exist in the interval $[0,\i)$.
\item(ii) If $d=4$ and Clebsch variables exist for $t=0$ outside some ball 
$B_r=\{x:|x|\le r\}$, then they exist globally in the interval $[0,\i)$. 
\end{itemize}
\end{proposition}

{\bf Proof. (i)} This is an immediate consequence of \eqref{eq1.12}: If $\Psi_0^*\bar\o=\o_0=d\a_0$, then
\[
\left(Q_t\circ \Psi_0^{-1}\right)^*d\a_t=\Psi_0^{-1*}Q_t^*d\a_t=\Psi_0^{-1*}d\a_0=\bar\o,
\]
which means that we can choose $\Psi_t=\Psi_0\circ A_t$ for the Clebsch change of variables.

{\bf (ii)} This is a consequence of a deep theorem of Gromov [Gr]: When $d=4$,
 a symplectic form is isomorphic to standard form $\bar\o$, if this is the case outside a ball $B_r$.
\qed

Observe that Euler Equation can be rewritten as
\begin{equation}\label{eq2.7}
\frac {d }{dt}\a_t+i_u(d\a_t)=-dH,
\end{equation}
where $H(q,t)=P(q,t)+|u(q,t)|^2/2$ is the Hamiltonian function. For a steady solution, $\a_t$ is independent of $t$ and we simply get
\[
i_u(d\a)=-dH.
\]
If $u$ is a symplectic steady solution of Euler Equation, then 
$i_u(d\a)=-dH$ means that $u$ is a {\it {Hamiltonian vector field}} with respect to the symplectic form $d\a$. Of course the associated the Hamiltonian function is $H$. Alternatively, we may write
\begin{equation}\label{eq2.7}
u=-\cC(u)^{-1}\nabla H.
\end{equation}

\begin{proposition}\label{prop2.2} Let $u$ be a steady symplectic solution to Euler Equation,
and let $c$ be a regular level set of $H(q,t)=P(q,t)+|u(q,t)|^2$ i.e. $\nabla H(q)\neq 0$ whenever $H(q)=c$.
Then the restriction of the form $\a$  to the submanifold $H=c$ is contact. In words, regular level sets of $H$ are contact submanifolds.
\end{proposition}

{\bf Proof.} By a standard fact in Symplectic Geometry (see for example [R]), the level set $H=c$ is contact if and only if we can find a {\it Liouville } vector field $X$ that is transversal to $M_c=\{H=c\}$. More precisely,
\[
\cL_X\ d\a=d\a,\ \ \ X(q)\notin T_qM_c,
\]
for every $q\in M_c$. Here $T_qM_c$ denotes the tangent fiber to $M_c$ at $q$. The first condition means that
$di_Xd\a=d\a$. This is  satisfied if $i_Xd\a=\a$. This really means that $\cC(u)X=u$ and as a result, we need to choose $X=\cC(u)^{-1}u$. It remains to show that $X$ is never tangent to $M_c$. For this, it suffices to check that $X\cdot \nabla H\neq 0$. Indeed, when $H=c$,
\[
X\cdot \nabla H=\cC(u)^{-1}u\cdot \nabla H=-u\cdot \cC(u)^{-1}\nabla H=|u|^2\neq 0,
\]
 by \eqref{eq2.8} because by assumption $\nabla H\neq 0$. We are done.
\qed

\noindent
{\bf{Example 2.1}} In this example we describe some simple solutions when the dimension is even.
We use polar coordinates to write $x_i=r_i\cos\th_i,$ $y_i=r_i\cos\th_i$, and let
$e_i$ (respectively $f_i$) denote the vector for which the $x_i$-th coordinate (respectively $y_i$-th coordinate)
is $1$ and any other coordinate is $0$. Set
\[
e_i(\th_i)=(\cos\th_i) e_i+(\sin\th_i) f_i,\ \ \ e'_i(\th_i)=(\sin\th_i) e_i-(\cos\th_i) f_i.
\]
We may write
\[
u=\sum_{i=1}^k\left(a^i e_i(\th_i)+b^ie'_i(\th_i)\right).
\]
The form $\a=u\cdot dx$ can be written as
\[
\a=\sum_{i=1}^k\left(a^i dr_i-r_ib^id\th_i\right)=:\sum_{i=1}^k\left(a^i dr_i-B^id\th_i\right).
\]
For a simple solution, let us assume that all $a^i$s and $b^i$s depend on $r=(r_1,\dots,r_d)$ only. 
We then have 
\[
d\a=\sum_{i<j}(a^i_{r_j}-a^j_{r_i})\ dr_i\wedge dr_j-\sum_{i,j}r_i^{-1}B^i_{r_j}\ dr_j\wedge (r_id\th_i).
\]
Now $u$ solves Euler Equation if the vector fields $a=(a^1,\dots,a^d)$ and $b=(b^1,\dots,b^d)$ satisfy
\begin{align}
&a_t+\cC(a)a-E(b)^*b+\nabla_r H=0,\nonumber\\
&b_t+E(b)a=0,\label{eq2.9}\\
&\sum_{i=1}^d (r_ia^i)_{r_i}/r_i=0,\nonumber
\end{align}
for some scalar function $H(r)$. Here $E(b)$ denotes 
a $d\x d$ matrix with entries $E_{ij}=r_i^{-1}B^i_{r_j}$. Note that if $\bar b=\sum_j (b^j)^2/2$, then
\[
E(b)^*b=[r_i^{-1}(b^i)^2]_i+\nabla \bar b=:\hat b+\nabla \bar b.
\]
Hence, by changing $H$ to $H'=H-\bar b$, we may rewrite the first equation in \eqref{eq2.9} with
\begin{equation}\label{eq2.10}
a_t+\cC(a)a-\hat b+\nabla_r H'=0.
\end{equation}
When $u$ is a steady solution, the first two equations in \eqref{eq2.9} simplifies to
\begin{equation}\label{eq2.11}
\cC(a)a-E(b)^*b+\nabla_r H=0,\ \ \ E(b)a=0.
\end{equation}
We can readily show that $u$ is a symplectic solution if and only if the matrix $E(b)$ is invertible.
Moreover, by taking the dot product of both sides of the first equation in \eqref{eq2.11},
and using the second equation we learn
\begin{equation}\label{eq2.12}
a\cdot \nabla H=0.
\end{equation}
Also, the equation $E(b)a=0$ really means
\begin{equation}\label{eq2.13}
a\cdot \nabla B^i=0 \ \ {\text{for }}i=1,\dots,d.
\end{equation}
When $d=4$ and $u$ is independent of time, it is straight forward to solve \eqref{eq2.9}: From the last equation in \eqref{eq2.9} we learn that there exists a function $\psi(r_1,r_2)$ such that 
\[
a^1=\psi_{r_2}/(r_1r_2),\ \ \ a^2=-\psi_{r_1}/(r_1r_2).
\]
From this, \eqref{eq2.12} and \eqref{eq2.13} we learn that $\nabla H$, $\nabla B^1$, $\nabla B^2$
 and $\nabla \psi$ are all parallel. So we may write
$$
H=\mu(\psi), \ \ B^1=\mu_1(\psi),\ \ B^2=\mu_2(\psi),
$$
 for some $C^1$ functions $\mu, mu_1,\mu_2:\bR\to\bR$.  Finally we go back to the first equation in
\eqref{eq2.11} to write
\[
a^2(a^2_{r_1}-a^1_{r_2})+\frac {B^1B^1_{r_1}}{r_1^2}+\frac {B^2B^2_{r_1}}{r_2^2}+H_{r_1}=0.
\]
Expressing this equation in terms of $\psi$  yields the elliptic PDE
\[
r_1r_2\left[\left(\frac {\psi_{r_1}}{r_1r_2}\right)_{r_1}+\left(\frac {\psi_{r_2}}{r_1r_2}\right)_{r_2}\right]
=\left(\frac {\mu'_1}{r_1^2}+\frac {\mu'_2}{r_2^2}-\mu'\right)(\psi).
\]
This equation may be compared to the Bragg-Hawthorne Equation that is solved to obtain axi-symmetric steady solutions in dimension three.
\qed

We now turn to the odd dimensions. assume that $d=2k+1$ for $k\in \bN$.
We write $(q_1,\dots,q_n)=(x_1,y_1,\dots,x_k,y_k,z)$ and when $k=1$
we simply write $(q_1,q_2,q_3)=(x,y,z)$. In this case, the standard 
contact form is $\bar\l=\sum_{i=1}^ky_idx_i+dz$. Again, locally all contact forms are isomorphic
to $\bar\l$.

\noindent
{\bf{Definition 2.3}}
Let $u$ be a  solution of Euler Equation. We say that {\it Clebsch} variables exist for $u$
in the interval $[0,T]$, if we can find $C^1$ functions
$$X_1,\dots,X_k,Y_1,\dots,Y_k:\bR^d\times [0,T]\to\bR,\ \ \ f,Z:\bR^d\times [0,T]\to\bR$$
such that $\Psi_t=(X_1,Y_1,\dots,X_k,Y_k,Z)(\cdot,t)$ is a diffeomorphism, $f>0$, and 
$$
(fu)(x,t)=\left(\sum_{i=1}^k Y_i\cdot \nabla X_i\right)(x,t)+\nabla Z(x,t),
$$
for every $t\in[0.T]$. 
Alternatively, we may write $f\a_t=\Psi_t^*\bar \l$.
\qed

As we recalled in the proof of Proposition 2.1(ii), if $d=4$ and 
a symplectic form is isomorphic to the standard form at
infinity, then the isomorphism can be extended to the whole $R^d$. This is no longer true when $d=3$; in fact there is countable collection of non-isomorphic forms $\l^n$ in $\bR^3$ such that each $\l^n$ is isomorphic to
$\bar \l$ at infinity but not globally. 
A fundamental result of Eliashberg gives a complete classification of contact forms. According to Eliashberg's Theorem [El], any contact form in $\bR^3$ is conformally isomorphic to one of the following forms
\begin{itemize}
\item (i) The standard form $\bar\l$.
\item (ii) The form $\hat \l=\frac {\sin r}{2r}(x_1  dx_2-x_2 dx_1)+{\cos r} \ dx_3$,
 where $r^2=x_1^2+x_2^2$.
\item (iii) A countable collection of pairwise non-isomorphic forms $\{\l^n:n\in\bZ\}$, where each $\l^n$ is isomorphic to $\bar \l$ outside the ball $B_1$ but not globally in $\bR^d$.
\end{itemize}
 
The above classification is related to the important notion of {\em{overtwisted}} contact forms. In fact $\hat \l$
is globally overtwisted whereas $\l^n$ are overtwisted only in a neighborhood of the origin.
(We refer to [El] or [Ge] for the definition of overtwisted forms).  

\noindent
{\bf{Example 2.2}} When $d=3$, we may use cylindrical coordinates $x_1=r\cos\th,x_2=r\sin\th$ to write
$u=ae(\th)+be'(\th)+ce_3,$
where
\[
e(\th)=r(\cos\th,\sin\th,0), \ \ e'(\th)=r(\sin\th,-\cos\th,0),\ \ e_3=(0,0,1).
\]
A solution is called axisymmetric if $a,b,$ and $c$ do not depend on $\th$. It turns out that any
\[
u=b(r)e(\th)+c(r)e_3,\ \ \a=u\cdot dx=c(r)\ dz-rb(r)\ d\th,
\]
is a steady solution to Euler equation. Such a solution is contact if 
\[
u\cdot\xi=r^{-1}(c(r)B'(r)-c'(r)B(r))\neq 0,
\] 
where $B(r)=rb(r)$. For example, if $b(r)=r,c(r)=1$, then we get
\[
\a=r^2\ d\th+dz=xdy-ydx+dz,
\]
is isomorphic to $\bar\l$. On the other hand, choosing $b(r)=r\sin r, c(r)=\cos r$ would yield exactly 
$\hat\l$.
\qed

\section{Symplectic Diffusions}
\label{sec3}
 \setcounter{equation}{0}

We study stochastic flows associated with diffusions. More precisely consider SDE
\begin{equation}\label{eq3.1}
dx(t)=V_0(x(t),t)dt+\sum_{i=1}^kV_i(x(t),t)\c dW^i(t),
\end{equation} 
where $(W^i:i=1,\dots,k)$ are standard one dimensional Brownian motions 
on some filtered probability space $(\O,\{\cF_t\},\bP)$, and $V_0,\dots,V_k$ are $C^r$--
vector fields in $\bR^n$. Here we are using Stratonovich stochastic differentials for the second term on the right-hand of \eqref{eq3.1} and a solution to the SDE \eqref{eq3.1} is a diffusion with the infinitesimal generator
\[
L=V_0\cdot\nabla+\frac 12 \sum_{i=1}^k(V_i\cdot\nabla)^2,
\]
 or in short $L=V_0+\frac 12\sum_{i=1}^kV_i^2$, where we have
simply written
$V$ for the $V$-directional derivative operator $V\cdot \nabla$.
We assume that the random flow $\phi_{s,t}$ of \eqref{eq3.1} is well defined almost surely. 
More precisely for $\bP-$almost all realization of $\o$, we have a flow $\{ \phi_{s,t}(\cdot,\o):0\le s\le t\}$
where $\phi_{s,t}(\cdot,\o):\bR^n\to\bR^n$ is a $C^{r-1}$ diffeomorphism and $\phi_{s,t}(a,\o)=:x(t)$ is a solution
of \eqref{eq3.1} subject to the initial condition $x(s)=a$. (We also write $\phi_t$ for $\phi_{0,t}$.)
For example a uniform bound on the $C^r$-norm of the coefficients
$V_0,\dots,V_k$ would guarantee the existence of such a stochastic 
flow provided that $r\ge 2$. We also remark that we can formally differentiate \eqref{eq3.1} with respect to the initial condition and derive a SDE for ${\L}_{s,t}(x)={\L}_t(x):=D_x\phi_{s,t}(x)$:
\begin{equation}\label{eq3.2}
d\L_t(x)=D_xV_0(\phi_{s,t}(x),t){\L}_t(x)dt+\sum_{i=1}^kD_xV_i(\phi_{s,t}(x),t)\ {\L}_t(x)\c dW^i(t).
\end{equation} 
Given a differential $\ell$-form $\a(x;v_1,\dots,v_\ell)$, we define
\[
\left(\phi_{s,t}^*\a\right)(x;v_1,\dots,v_\ell)=\a(\phi_{s,t}(x);{\L}_{s,t}(x)v_1,\dots,{\L}_{s,t}(x)v_\ell).
\]
Given a vector field $V$, we write $\cL_V$ for the Lie derivative in the direction $V$. More precisely, for every differential 
form $\a$,
\begin{equation}\label{eq3.3}
\cL_V\a=(\hat d\c i_V+i_V\c \hat d)\a,
\end{equation}   
where $\hat d$ and $i_V$ denote the exterior derivative and $V-$contraction operator respectively. 
(To avoid a confusion between the stochastic differential and exterior derivative, we are using a hat for the latter.)
We are now ready to state 
a formula that is the stochastic analog of Cartan's formula and it is a rather straight forward consequence of \eqref{eq3.2}. We refer to Kunita [K2] for a proof. 
\begin{proposition}
\label{prop3.1} Set ${\bf {V}}=(V_0,V_1,\dots,V_m)$ and
\[
\cA_{\bf V}=\cL_{V_0}+\frac 12\sum_{i=1}^k\cL_{V_i}^2.
\]
We also $\eta_t$ for $\phi_{s,t}^*\eta$ for any form $\eta$. We have
 \begin{align}
d\a_t&=\left(\cL_{V_0}\a\right)_t\ dt+\sum_{i=1}^k\left(\cL_{V_i}\a\right)_t \c dW^i(t)\label{eq3.4}\\
&=\left(\cA_{\bf V}\a\right)_t \ dt+\sum_{i=1}^k\left(\cL_{V_i}\a\right)_t \ dW^i(t)\nonumber.
\end{align}
\end{proposition}

\noindent
{\bf Example 3.1}
\begin{itemize}
\item (i) If $\a=f$ is a $0$-form, then $\cA_{\bf V}f=Lf$ is simply the infinitesimal generator of the underlying diffusion.
\item (ii) If $\a=\rho\ dx_1\wedge\dots\wedge dx_n$, is a volume form, then $\cA_{\bf V}\a=(L^*\rho)\ dx_1\wedge\dots\wedge dx_n$, where $L^*$ is the adjoint of the operator $L$.
\item (iii) If $W=(W^1,\dots,W^n)$ is a $n$-dimensional standard Brownian motion and 
 \[
dx=V_0(q,t)dt+dW,
\]
then for  a volume form $\a=\rho\ dx_1\wedge\dots\wedge dx_n$, we write
 $\a_t=\rho^t\ dx_1\wedge\dots\wedge dx_n$, and  \eqref{eq3.4} becomes
\[
d\rho^t=L^*\rho^t\ dt+\nabla\rho^t\cdot dW.
\]
In particular, when $\nabla\cdot V_0=0$ and $\rho=\rho^0=1$, then $\rho^t=1$ is a solution. In other words,
the standard volume $dx_1\wedge\dots\wedge dx_n$ is preserved for such a diffusion if the drift $V_0$ is divergence free.
\end{itemize}
\qed

 We now make two definitions:
\begin{definition}\label{def3.1}
Let $\a$ be a symplectic form.
\begin{itemize}
\item (i) We say that the diffusion \eqref{eq3.1} is (strongly) $\a$-symplectic if its flow is symplectic with respect $\a$, 
almost surely. That is $\phi_{t}^*\a=\a$, a.s.
\item(ii) We say that the diffusion \eqref{eq3.1} is weakly symplectic if  $\a_t:=\phi_{t}^*\a$, is a martingale.
\end{itemize}
\end{definition}

Using Proposition 3.1 it is not hard to deduce

\begin{proposition}\label{prop3.2}
\begin{itemize}
\item (i) The diffusion \eqref{eq3.1} is (strongly) $\a$-symplectic if and only if the vector fields $V_0,V_1,\dots,V_k$ are $\a$-Hamiltonian, i.e. $\cL_{V_0}\a=\cL_{V_1}\a=\dots=\cL_{V_k}\a=0$.
\item(ii) The diffusion \eqref{eq3.1} is weakly $\a$-symplectic if and only if $\cA_{\bf V}\a=0$.
\end{itemize}
\end{proposition}

We discuss two systematic ways of producing weakly symplectic diffusions. 

{\bf {Recipe (i)}} Given a symplectic form $\a$,
we write $X_H=X_H^\a$ for the Hamiltonian vector field associated with the Hamiltonian function $H$. Note that by non-degeneracy of $\a$, there exists a unique vector field $X=\cX^\a(\nu)$ such that $i_X \a=\nu$ for every $1$-form 
$\nu$ and $X_H=-\cX^\a(dH)$. In the following proposition, we show that given $V_1,V_2,\dots,V_k$, we can always find
a unique $\hat V_0$ such that the diffusion associated with ${\bf V}
=(X_H+\hat V_0,V_1,\dots,V_k)$ is weakly $\a$-symplectic.

\begin{proposition}\label{prop3.3}
 The diffusion \eqref{eq3.1} is weakly $\a$-symplectic if and only if there exists a Hamiltonian function $H$, such that 
\begin{equation}\label{eq3.5}
V_0=X_H-\frac 12\sum_{j=1}^k\cX^\a\left(i_{V_j}\ \hat d\ i_{V_j}\a\right).
\end{equation}
\end{proposition}

{\em{Proof.}} By definition,
\[
\cA_{{\bf V}}\a=\hat d\left[i_{V_0}\a+\frac 12\sum_{j=1}^k\left(i_{V_j}\ \hat d\ i_{V_j}\a\right)\right].
\]
Hence  $\cA_{{\bf V}}\a=0$ means that for some function $H$, 
\[
i_{V_0}\a+\frac 12\sum_{j=1}^k\left(i_{V_j}\ \hat d\ i_{V_j}\a\right)=-dH.
\]
From this we can readily deduce \eqref{eq3.5}.
\qed

{\bf {Recipe (ii)}} We now give a useful recipe for constructing $\bar\o$-diffusions where $\bar\o$ is the standard symplectic form and $n=2d$. 

\begin{proposition}\label{prop3.3}
 Given a Hamiltonian function $H$, consider a diffusion $x(t)=(q(t),p(t))$ that solves
\begin{align}
dq&=H_p(q,p)\ dt+\sum_{j=1}^kA_j(x,t)\ dW^j,\label{eq3.6}\\
dp&=-H_q(q,p)\ dt+\sum_{j=1}^kB_j(x,t)\ dW^j,\nonumber
\end{align}
with $A_j=(A_j^1,\dots,A_j^n)$, and $B_j=(B_j^1,\dots,B_j^n)$.
Then $x(t)$ is weakly $\bar\o$-symplectic if and only if $Z_1=(Z_1^i:i=1,\dots,d)=0$ and $Z_2=(Z_2^i:i=1,\dots,d)=0$,
where
\begin{align}
Z_1^i&=\sum_{r,j}\left(\frac{\partial A^r_j}{\partial q_i}B_j^r
-\frac{\partial B^r_j}{\partial q_i} A_j^r\right),\label{eq3.7}\\
Z_2^i&=\sum_{r,j}\left(\frac{\partial A^r_j}{\partial p_i}B^r_j-\frac{\partial B^r_j}{\partial p_i}A^r_j\right). \nonumber
\end{align}

\end{proposition}
{\em{Proof.}}
The Stratonovich differential is related to It\^o differential by 
\[
a\circ dW=a\ dW+\frac 12 [da,dW].
\]
As a result, the diffusion $x(t)$ satisfies \eqref{eq3.1} for $V_j=\bmatrix A_j\\B_j\endbmatrix$, $j=1,\dots,k$, and
$V_0=J\nabla H-\frac 12\hat V_0$ with $\hat V_0=\bmatrix A_0\\B_0\endbmatrix$, where
\begin{align}
A_0^i&=\sum_{r,j}\left(\frac {\partial A_j^i}{\partial q_r }A^r_j+
\frac {\partial A^i_j}{\partial p_r }B_j^r\right) ,\label{eq2.8}\\
B_0^i&=\sum_{r,j}\left(\frac {\partial B_j^i}{\partial q_r }A^r_j+
\frac {\partial B^i_j}{\partial p_r }B_j^r\right) .\nonumber
\end{align}
We need to show that \eqref{eq3.5} is satisfied if and only if $Z_1=Z_2=0$. For this, 
let us write $\b( F)$ for the $1$-form $F\cdot dx$ and observe 
\[
i_V\bar\o=\b(JV),\ \ \ \hat d\b(F)(v,w)=\cC(F)v\cdot w,\ \ \ \cX^{\bar\o}(\b(F))=-J F,
\]
where $\cC(F)=DF-(DF)^*$ with $DF$ denoting the matrix of the partial derivatives of $F$ with respect to $x$.
From this we deduce
 \[
\sum_{j=1}^k\cX^{\bar\o}\left(i_{V_j}\ \hat d\ i_{V_j}\bar\o\right)=-\sum_{j=1}^kJ\cC(JV_j)V_j.
\]
On account of this formula and Proposition 3.3, it remains to verify 
\begin{equation}\label{eq3.9}
\hat V_0=-\sum_{j=1}^kJ\cC(JV_j)V_j.
\end{equation}

A straight forward calculation yields
\[
\cC(JV_j)=\bmatrix
X_{11}^j & X_{12}^j \\
X_{21}^j & X_{22}^j\endbmatrix
\]
where
\begin{align*}
X_{11}^j&=\left[\frac{\partial B^i_j}{\partial q_r}-\frac{\partial B^r_j}{\partial q_i}\right]_{i,r=1}^n,
\ \ \ X_{12}^j=\left[\frac{\partial B^i_j}{\partial p_r}+\frac{\partial A^r_j}{\partial q_i}\right]_{i,r=1}^n,\\
 X_{21}^j&=\left[-\frac{\partial A^i_j}{\partial q_r}-\frac{\partial B^r_j}{\partial p_i}\right]_{i,r=1}^n,
\ \ \ X_{22}^j=\left[-\frac{\partial A^i_j}{\partial p_r}+\frac{\partial A^r_j}{\partial p_i}\right]_{i,r=1}^n.
\end{align*}
From this we deduce 
\[
J\cC(JV_j)V_j=\bmatrix Y^j_1\\Y^j_2\endbmatrix,
\]
where
\begin{align*}
Y_1^j&=\left[-\sum_r\left(\frac{\partial A^i_j}{\partial q_r}+\frac{\partial B^r_j}{\partial p_i}\right)A^r_j
+\sum_r\left(\frac{\partial A^r_j}{\partial p_i}-\frac{\partial A^i_j}{\partial p_r}\right)B^r_j\right]_{i=1}^n,\\
Y_2^j&=\left[\sum_r\left(\frac{\partial B^r_j}{\partial q_i}-\frac{\partial B^i_j}{\partial q_r} \right)A_j^r
-\sum_r\left(\frac{\partial B^i_j}{\partial p_r}+\frac{\partial A^r_j}{\partial q_i}\right)B_j^r\right]_{i=1}^n.
\end{align*}
Summing these expressions over $j$ yields
\[
-\sum_jJ\cC(JV_j)V_j=\hat V_0-\bmatrix Z_2\\-Z_1\endbmatrix=V_0-J\bmatrix Z_1\\Z_2\endbmatrix,
\]
where $Z_1$ and $Z_2$ are defined by \eqref{eq3.7}.
From this we learn that \eqref{eq3.9} is valid if and if $Z_1=Z_2=0$. This completes the proof.
\qed

An immediate consequence of Proposition 3.4 is Corollary 3.1.
\begin{corollary}\label{cor3.1}
Let $x(t)=(q(t),p(t))$ be a diffusion satisfying
\begin{align}
dq&=H_p(q,p)\ dt+\sqrt{2\nu}\ dW ,\label{eq3.10}\\
dp&=-H_q(q,p)\ dt +\sqrt{2\nu}\G(q,t) dW.\nonumber
\end{align}
where $\G$ is a continuously differentiable $d\x d$-matrix valued function and $W=(W^1,\dots,W^d)$ is a standard Brownian motion in $\bR^d$.
The process $x(t)$ is weakly $\bar\o$-symplectic if and only if the trace of $\G$ is independent of $q$.
\end{corollary}
{\em{Proof.}} Observe that $x(t)$ satisfies \eqref{eq3.6} for $A=I_d$ and $B$ that is independent of $p$.
From this we deduce that $Z_2=0$ and $Z_1=-\sqrt{2\nu}\nabla_q ( tr \G)$. We are done.
\qed
\section{Martingale Circulation}
\label{sec4}
 \setcounter{equation}{0}

{\em{Proof of Theorem 1.1.}} {\bf{Step 1.}} 
As in Section 1, we write $D$ and $\nabla$ for $q$-differentiation.
  For $x$-differentiation however, we write $D_x$ and $\nabla_x$ instead.
Let us write $x'(t)=(q'(t),p'(t))$ for a diffusion that satisfies
\begin{align}
dq'(t)&=p'(t)\ dt+\sqrt{2\nu}\ d\bar W\nonumber\\
dp'(t)&=-\nabla P(q'(t),t) \ dt+\sqrt{2\nu}\ Dw(q'(t),t)\ d\bar W.\label{eq4.1}
\end{align}
for a time dependent $C^1$ vector field $w$ in $\bR^d$ and a standard Brownian motion $\bar W$. 
The flow of this diffusion is denoted by $\phi_t$. 
We then apply Corollary 3.1 for $H(q,p,t)=\frac 12 |p|^2+P(q,t)$ and
$\G=Dw$, to assert that the 
diffusion $x'$ is weakly $\bar\o$-symplectic if $\nabla \cdot w=0$.
Let us now assume that $w$  satisfies the backward Navier-Stokes equation
\begin{equation}\label{eq4.2}
w_t+(Dw)w+\nabla P+\nu\D w=0,\ \ \  \nabla\cdot w=0.
\end{equation} 
We observe that if the process $q'(t)$ is a diffusion satisfying
\begin{equation}\label{eq4.3}
dq'(t)=w(q'(t),t)\ dt+\sqrt{2\nu}\ d\bar W.
\end{equation}
and $p'(t)=w(q'(t),t)$, then by Ito's formula,
\begin{align}
dp'(t)&=\left[w_t+(Dw)w+\nu \D w\right](q'(t),t) \ dt+\sqrt{2\nu}Dw(q'(t),t)\ d\bar W\nonumber\\
&=-\nabla P(q'(t),t) \ dt+\sqrt{2\nu}Dw(q'(t),t)\ d\bar W.\label{eq4.4}
\end{align}
This means that if $\bar Q_t$ denotes the flow of the SDE \eqref{eq4.3}, then
\begin{align}\label{eq4.5}
\phi_t(q,w(q,0))&=(\bar Q_t(q),w(\bar Q_t(q),t)),\\
D_x\phi_t{(q,w(q,0))}\bmatrix a\\Dw(q,0)a\endbmatrix&=\bmatrix (D\bar Q_t(q))a\\(Dw(\bar Q_t(q),t))
(D\bar Q_t(q))a\endbmatrix.\nonumber
\end{align}
By the conclusion of Corollary 3.1, the process 
\[
\hat M_t(x;v_1,v_2)=\left[J(D\phi_t(x))v_1\right]\cdot \left[(D\phi_t(x))v_2\right]
\]
is a $2$-form valued martingale.  This means that for any 
surface $\g:D\to\bR^d\x \bR^d$, the process
\[
\hat M_t(\g)=\iint_D \hat M_t(\g;\g_{\th_1},\g_{\th_2})\ d\th_1 d\th_2,
\]
is a martingale. We consider a surface that lies on the graph of $w(\cdot,0)$. That is, 
$$\g(\th_1,\th_2)=(\tau(\th_1,\th_2), w(\tau(\th_1,\th_2),0)),$$
for a surface $\tau:D\to\bR^d$. We now use \eqref{eq4.5} to assert that 
\[
\bar M_t(\tau):=\hat M_t(\g)=\iint_D \bar  M_t(\tau;\tau_{\th_1},\tau_{\th_2})\ d\th_1 d\th_2,
\]
where
\begin{align*}
\bar M_t(q;a_1,a_2)&=J\bmatrix (D\bar Q_t(q))a_1\\(Dw(\bar Q_t(q),t))(D\bar Q_t(q))a_1\endbmatrix\cdot
 \bmatrix (D\bar Q_t(q))a_2\\(Dw(\bar Q_t(q),t))(D\bar Q_t(q))a_2\endbmatrix\\
&=\left[(Dw-(Dw)^*)( \bar Q_t(q),t)\right](D\bar Q_t(q))a_1\cdot (D\bar Q_t(q))a_2\\
&=\bar Q_t^*\hat d\bar \a_t(q;a_1,a_2),
\end{align*}
where $\bar\a_t=w(q,t)\cdot dq$. In summary, $\bar M_t=\bar Q_t^*\hat d\bar \a_t$ is a martingale.

When $d=3$,
\begin{align*}
\bar Q_t^*\hat d\bar \a_t(q;a_1,a_2)
 &=\left[(\eta^t\circ \bar Q_t(q))\x (D\bar Q_t(q))a_1\right]\cdot (D\bar Q_t(q))a_2\\
&=\left[\eta^t\circ \bar Q_t(q),(D\bar Q_t(q))a_1,(D\bar Q_t(q))a_2\right],
\end{align*}
where $\eta^t(\cdot)=\nabla\x w(\cdot,t)$ and $[a,b,c]$ is the determinant of a matrix with column 
vectors $a,b$ and $c$. Since $w$ is divergence-free, the flow $\bar Q_t$ is volume preserving
(see Example 3.1(iii)). Hence
\[
\bar M_t(q;a_1,a_2)=[(D\bar A_t\circ \bar Q_t(q))\eta^t\circ \bar Q_t(q),a_1,a_2],
\]
where $\bar A_t=\bar Q_t^{-1}$. Since $\hat M_t$ is a martingale, we deduce that the process
\[
\tilde M_t(q)=(D\bar A^t\circ \bar Q_t(q))(\eta^t\circ \bar Q_t(q)),
\]
is a martingale.

{\bf{Step 2.}} Suppose that now $u$ is a solution to the forward Navier-Stokes equation \eqref{eq1.15}
 and recall that when $d=3$, we write
$\xi=\nabla\x u$. We set $w(q,t)=-u(q,T-t)$ for $t\in[0,T]$. Then $w$ satisfies \eqref{eq3.2} in the interval
$t\in[0,T]$. Recall that 
$q(t)$ is the solution of SDE \eqref{eq1.16} with the flow $Q_t$. We choose $\bar W(t)=W(T-t)-W(T)$ 
in the equation
\eqref{eq4.3}. According to a theorem of Kunita (see Theorem 13.15 in page 139 of [RW]
and [K1]), the flows $Q$ and $\bar Q$ are related by the formula
\[
\bar Q_t=Q_{T-t}\circ Q_T^{-1}=B_t.
\]
Observe that $\bar \a_t=-\a_{T-t}$ and 
\[
\bar M_t=\bar Q_t^*\hat d\bar\a_t=-B_t^*\hat d\a_{T-t}=-\b_t.
\]
Hence $(\b_t: t\in[0,T])$ is a  martingale because $\bar M_t$ is a martingale by Step 1.
 Also, when $d=3$,
\[
\tilde M_t=\big((D\bar A_t)\eta^t\big)\circ \bar Q_t(q))=-
\left(\left(DB_t^{-1}\right)\xi^{T-t}\right) \circ B_t.
\]
 This completes the proof of Part $(i)$.

{\bf{Step 3.}} The process $x'(t)$ is a diffusion of the form \eqref{eq4.1} with $k=d$ and 
\[
V_i(x',t)=V_i(q,t)=\bmatrix e_i\\w_{q_i}\endbmatrix,
\]
 where $e_i=[\d_i^j]_{j=1}^d$ is the unit vector in the $i$-th direction.
A straight forward calculation yields that for the standard symplectic form $\bar\o=\sum_j dp_j\wedge dq_j$,
\begin{align*}
i_{V_i}\bar\o&=w_{q_i}\cdot dq-dp_i=:\g^i-dp_i,\\ 
\cL_{V_i}\ \bar\o&=\hat d\g^i=\sum_{j,k}w^k_{q_iq_j}\ dq_j\wedge dq_k,\\
\zeta^i(v_1,v_2)&=\hat d\g^i(v_1,v_2)=\cC(w_{q_i})v_1\cdot v_2,
\end{align*}
where $w=(w^1,\dots,w^d)$.  From this and \eqref{eq3.4} we deduce if 
\[
z_t=\int_\Theta\ (\hat d\bar \a)_t,\ \ \ y_i(t)=\int_\Theta\ \zeta^i_t,
\]
 then 
\[
dz_t=\sum_{i=1}^d y_i(s)\ dW^i(t),
\]
because by Step 1, we know that $\cA_{\bf V}\o=0$.
From this, we readily deduce that the quadratic variation of the process
$z_t$ is given by
\[
\int_0^t \sum_i y_i(s)^2\ ds.
\]
We now reverse time as in Step 2 to complete the proof of Part $(ii)$.
Part $(iii)$ is an immediate consequence of the identity
\[
\bE z_T^2=\bE z_0^2+\bE\int_0^T\sum_i y_i(s)^2\ ds.
\]
\qed

\section{Contact Diffusions}
\label{sec5}
 \setcounter{equation}{0}

Recall that contact forms are certain $1$-forms that are non-degenerate in some rather strong sense.
 To explain this, recall that when $\a$ is a contact form in dimension $n=2d+1$, then 
the set $l_x=\{v:d\a_x(v,w)=0 \text{ for all }w\in T_xM\}$ is a line. Also, if we define the kernel of $\a$
by
\begin{equation*}
\eta_x^\a=\eta_x=\{v: \a_x(v)=0\},
\end{equation*}
then  the contact condition really means that $l_x$ and $\eta_x$ give a decomposition of $\bR^{n}$ that depends solely on $\a$:
\begin{equation}\label{eq5.1}
\bR^{n}=\eta_x\oplus l_x.
\end{equation}
We also define the Reeb vector field $R(x)=R^\a(x)$ to be the unique vector such that 
\begin{equation*}
R(x)\in l_x,\ \ \  \a_x(R(x))=1.
\end{equation*}
 
The role of Hamiltonian vector fields in the contact geometry are played by contact vector field.

\noindent
{\bf{Definition 5.1}} A vector field $X$ is called an $\a$-contact vector field if $\cL_X\a=f\a$ for some scalar-valued continuous function $f$.
\qed

It is known that for a
 given a ``Hamiltonian'' $H:M\to\bR$, there exists a unique contact $\a$-vector field $X_H=X_{H,\a}$ 
such that $i_{X_H}\a=\a({X_H})=H$. The function $f$ can be expressed in terms of $H$ with the aid of the Reeb's vector field $R=R^\a$; indeed, $f=dH(R^\a)$, and as a result,
\[
\cL_{X_H}\a=dH(R^\a)\a.
\]
In our Euclidean setting, we consider a form
 $\a=u\cdot dx$ for a vector field $u$ and 
\[
\b(v_1,v_2):=d\a(v_1,v_2)=\cC(u)v_1\cdot v_2,
\]
where $\cC(u)=Du-(Du)^*$. (Recall that we are writing $A^*$ for the transpose of $A$.) 
Since $C^*=-C$, we have that
$\det \cC=(-1)^n \det \cC$. This implies that $C$ cannot be invertible if the dimension is odd.
 Hence the null space $l_x$ of $\cC(u)(x)$ is never trivial and our assumption $dim\ell_x=1$ really means
that  this null space has the smallest possible dimension. Now \eqref{eq5.1} simply means that $u(x)\cdot R(x)\neq 0$.  Of course $R$ is chosen so that $u(x)\cdot R(x)\equiv 1$. Writing $u^{\perp}$ and $R^{\perp}$ for the space of vectors perpendicular to $u$ and $R$ respectively, then $\eta=u^\perp$, and
we may define a matrix $\cC'(u)$ which is not exactly the inverse of $\cC(u)$ (because 
$\cC(u)$ is not invertible),  but it is specified uniquely by two requirements:
\begin{itemize} 
\item (i) $\cC'(u)$ restricted to $R^\perp$ is the inverse of $\cC(u):u^\perp\to R^\perp$.
\item (ii) $\cC'(u)R=0$. 
\end{itemize}
The contact vector field associated with $H$ is given by
\begin{equation*}
X_H=-\cC'(u) \nabla H+H R.
\end{equation*}
 
In particular, when $n=3$, the form $\a=u\cdot dx$
is contact if and only if $u\cdot \xi$ is never $0$, where $\xi=\nabla \x u$ is the curl (vorticity) of $u$.
In this case the Reeb vector field is given by $R=\xi/(u\cdot \xi)$, and
\begin{equation*}
\cL_Z u=\nabla(u\cdot X)+\xi\x Z,
\end{equation*} 
We also write $\bar u=u/\rho$. The contact vector field associated with $H$ is given by
\begin{equation*}
X_H=\bar u\x \nabla H+H R.
\end{equation*}

Let $x(t)$ be a diffusion satisfying \eqref{eq3.1} and assume that this diffusion has a random flow $\phi_t$.
Given a contact form $\a=\a_0$, set $\a_t=\phi_t^*\a_0$ as before. 

 \noindent
{\bf Definition 5.2.} 
\begin{itemize}
\item (i) We say that the diffusion \eqref{eq3.1} is strongly $\a$-contact, if for some scaler-valued semimartingale $Z_t$ of the form,
\begin{equation}\label{eq5.2}
dZ_t=g_0(x(t),t)\ dt+\sum_{i=1}^kg_i(x(t),t)\circ dW^i(t),
\end{equation}
we have
\[
d\a_t=\a_t\ dZ_t.
\]
\item (ii) We say that the diffusion \eqref{eq3.1} is weakly $\a$-contact, if there exists a continuous
scalar-valued function $f(x,t)$ such that
\[
M_t=\a_t-\int_0^tf(x(s),s)\a_s\ ds,
\]
is a martingale. 
\end{itemize}
\qed
  
We end this section with two proposition. 

\begin{proposition}\label{prop5.1} The following statements are equivalents:
\begin{itemize}
\item(i) The diffusion \eqref{eq3.1} is strongly $\a$-contact.
\item (ii) There exists a scaler-valued process $A_t$ of the form
\[
dA_t=h_0(x(t),t)\ dt+\sum_{i=1}^kh_i(x(t),t)\circ dW^i(t).
\]
 such that $\a_t=e^{A_t}\a.$ (Recall $\a_t=\phi_t^*\a$ with $\phi_t=\phi_{0,t}$ representing the flow of the diffusion \eqref{eq3.1}.)
\item(iii) The vector fields $V_0,\dots,V_k$ are $\a$-contact.
\end{itemize}
\end{proposition}

\begin{proposition}\label{prop5.2} The following statements are equivalents:
\begin{itemize}
\item(i) The diffusion \eqref{eq3.1} is weakly $\a$-contact.
\item(ii) For some scalar-valued function $f(x,t)$, we have $\cA_{\bf V}\a=f\a$.
\end{itemize}
\end{proposition}

The proof of Proposition 5.2 is omitted because it is an immediate consequence of \eqref{eq3.4} and the definition.

{\bf{Proof of Proposition 5.1.}} Suppose that the vector fields $V_0,\dots,V_k$ are $\a$-contact. Then there exist scalar-valued functions
$g_0(x,t),\dots,g_k(x,t)$ such that $\cL_{V_i}\a=g_i\a$. From this and Proposition 3.1 we learn that 
$d\a_t=\a_t \ dZ_t$ for
$Z_t$ as in \eqref{eq5.2}.
Hence $(iii)$ implies $(i)$. 

Now assume $(i)$ and set
\[
Y_t=\exp\left(-Z_t+\frac 12[Z]_t\right).
\]
 We have 
\begin{align*}
dY_t&=Y_t(- dZ_t+d[Z]_t),\\
d\left(Y_t \a_t\right)&=\a_tY_t(- dZ_t+d[Z]_t)+Y\ d\a_t+ d[Y,\a]_t\\
&= \a_tY_t\ d[Z]_t+ d[Y,\a]_t= \a_tY_t\ d[Z]_t-
 \a_tY_t\ d[Z]_t=0.
\end{align*}
Hence $Y_t\a_t=\a$ and we have $(ii)$ for $A_t=Z_t-\frac 12[Z]_t$.

We now assume $(ii)$. We certainly have
\begin{align*}
d\a_t&=\a e^{A_t}\left(dA_t+\frac 12d[A]_t\right)=\a_t\left(dA_t+\frac 12d[A]_t\right)\\
&=\a_t\left(g_0(x(t),t)\ dt+\sum_{i=1}^kg_i(x(t),t)\circ dW^i(t)\right),
\end{align*}
for $g_0=h_0+\frac 12(\sum_i h_i^2)$ and $g_i=h_i$ for $i=1,\dots,k$.
Comparing this to \eqref{eq3.4} yields  $\cL_{V_i}\a=g_i\a$ for $i=0,\dots,k$. Hence $(iii)$ is true and this completes the proof.
\qed


\begin{thebibliography}{HR}



\bibitem[CI]{CI}P. Constantin and G. Iyer,  A stochastic Lagrangian representation of the
three-dimensional incompressible Navier-Stokes equations. Commun. Pure
Appl. Math. {\bf LXI}, 330-345, (2008).
\bibitem[El]{El}Y. Eliashberg, Classification of contact structures on $\bR^3$. 
Internat. Math. Res. Notices , no. 3, 87–91, (1993).

\bibitem[Ey1]{E1}G. L. Eyink,  Turbulent diffusion of lines and circulations. Physics
Letters A {\bf 368}, 486-490, (2007).
\bibitem[Ey2]{Ey2}G. L. Eyink, A stochastic least-action principle for the incompressible
Navier-Stokes equation, submitted to Physica D. arxiv:0810.0817 [math-ph], (2009). 
\bibitem[Ge]{Ge} H. Geiges, An Introduction to Contact Topology.
 Cambridge university Press, Cambridge, 2008.

\bibitem[Gr]{Gr} M. Gromov, Pseudoholomorphic curves in symplectic manifolds. Invent. Math. {\bf 82},  
307–347, (1985).

\bibitem[HZ]{HZ} H. Hofer and E. Zehnder, Symplectic Invariants and Hamiltonian Dynamics. 
Birkhäuser Verlag, Basel, 2011. 

\bibitem[K1]{K1} H. Kunita, On backward stochastic differential equations. Stochastica, {\bf 6}, 293-313 (1982).

\bibitem[K2]{K2} H. Kunita, Stochastic Flows and Stochastic Differential Equations.
Cambridge university Press, Cambridge, 1990.

\bibitem [R]{R} F.  Rezakhanlou, Lectures on Symplectic Geometry,  www.math.berkeley.edu/rezakhan.

\bibitem[RW]{RW} L. C. G. Rogers and D. Williams,  Diffusions, Markov processes, and martingales. Vol. 2. It\^o calculus.  Cambridge University Press, Cambridge, 2000.


 
      

\end{thebibliography}
\end{document}